\documentclass[10pt,a4paper]{amsart}
\usepackage{amsfonts}
\usepackage{epsfig}
\usepackage{graphicx}
\usepackage{amsmath}
\usepackage{amssymb}
\usepackage{mathrsfs}
\usepackage{color}
\newcommand{\Rr}{\mathbb{R}}

\newcommand{\R}{\mathbb{R}}

\newcommand{\Z}{\mathbb{Z}}

\DeclareMathAlphabet{\mathpzc}{OT1}{pzc}{m}{it}

\newcounter{main}

\newtheorem{thm}{Theorem}[section]
\newtheorem{lem}[thm]{Lemma}
\newtheorem{pro}[thm]{Proposition}

\newtheorem{maintheorem}{Theorem}
\newtheorem{cor}[thm]{Corollary}
\newtheorem{con}[thm]{Conjecture}
\newtheorem{rk}[thm]{Remark}

\DeclareMathAlphabet{\mathpzc}{OT1}{pzc}{m}{it}
\DeclareSymbolFont{symbols}{OMS}{cmsy}{m}{n}
\usepackage[T1]{fontenc}

\title[Shadowing, expansiveness and specification for $C^1$-conservative systems]{Shadowing, expansiveness and specification for $C^1$-conservative systems}

\author[M. Bessa]{M\'{a}rio Bessa}

\address{Universidade da Beira Interior, Rua Marqu\^es d'\'Avila e Bolama,
  6201-001 Covilh\~a
Portugal.}
\email{bessa@fc.up.pt}

\author[M. Lee]{Manseob Lee}

\address{Department of Mathematics
Mokwon University, Daejeon 302-729, Korea.}
\email{lmsds@mokwon.ac.kr}

\author[X. Wen]{Xiao Wen}

\address{The School of Mathematics and Systems Sciences, Beihang University, Beijing 100191, P. R. China}
\email{wenxiao@buaa.edu.cn}

\begin{document}

\begin{abstract}
We prove that a $C^1$-generic volume-preserving dynamical system (diffeomorphism or flow) has the shadowing property or is expansive or has the weak specification property if and
only if it is Anosov. Finally, as in \cite{B2,F2}, we prove that the $C^1$-robustness, within the volume-preserving context, of the expansiveness property and the weak specification property, imply that the dynamical system (diffeomorphism or flow) is Anosov.
\end{abstract}

\maketitle

\noindent{ 2010 {\em Mathematical Subject Classification}: 37C20, 37C50. \\
{\em Keywords and phrases}: shadowing, expansiveness, specification, generic, Anosov,
vo\-lu\-me-preserving, star systems.}

\section{Introduction, basic definitions and statement of the results}

\subsection{Introduction}
One of the main goals in dynamics is the understanding of the orbits
structure for most maps and flows. Since it is an unattainable task understand \emph{all} dynamical systems people often try to understand large classes of systems. The natural decreasing hierarchy among large subsets of all systems is: 
\begin{itemize}
\item The \emph{richer}: an \emph{open and dense} subset,  
\item The \emph{in between}: a \emph{generic} (in the sense of second Baire's category) and 
\item The \emph{poorer}: a \emph{dense} subset.  
\end{itemize}
Fixed a certain topology in the class of a Baire space of dynamical systems it is very important to detect when a dense property holds generically with respect to that topology. The reason is simple; a residual subset is not only ``thicker'' than a dense one, but countable intersections of residual subsets is still residual and this property is no longer true for dense subsets.

The notion of shadowing in dynamics is intuited by the numerical computational idea of evaluate discrepancies between exact and approximate solutions along orbits and to perceive the influence of the errors that we commit and allow on each iterate. In few words the \emph{shadowing property} aims to obtain shadowing of approximate trajectories in a given dynamical system by true orbits of the system. 

The concept of \emph{expansiveness} means, in rough terms, that if two points stay near for forward and backward iterates, then they must be equal. In some sense, expansive systems can be considered chaotic since they exhibit sensitivity to the initial conditions. 

Informal speaking the notion of \emph{specification}  of a dynamical system means that it is possible to shadow two distinct ``pieces of orbits'', which are sufficiently apart in time, by a single orbit. 

All previous notions are central in hyperbolic dynamics and reveal to be very useful in the confluence of dynamics and differentiable ergodic theory (see ~\cite{Bow,DGS}).

It is well known that the Anosov diffeomorphisms satisfy the three aforementioned properties. 
With respect to stability we recall, for instance, that in \cite{R1}
Robinson showed that a structurally stable diffeomorphism has
the shadowing property, and in \cite{S} Sakai proved that the
$C^1$-interior of the set of diffeomorphisms with the shadowing
property coincides with the set of structurally stable
diffeomorphisms. We refer also recent results by one of the authors (\cite{WGW}) and also ~\cite{P1,LS} for flows. Inward  the dissipative setting it was proved that  $C^1$-stability of shadowing and also of expansiveness guarantee that the systems is hyperbolic (\cite{M0,S2,MSS,ST}).

In the present paper we intend to study the class of volume-preserving (or conservative) diffeomorphisms and flows endowed with the $C^1$-Whitney topology. Namely, we consider that the manifold has dimension greater or equal than three for the discrete case, and we consider that the manifold has dimension greater or equal than four for the flow case.
It is already known that, if a volume-preserving diffeomorphism $f$ has the shadowing property and, moreover, any volume-preserving diffeomorphism arbitrarily $C^1$-close to $f$ displays also the shadowing property (see \cite{B2}), then $f$ must be Anosov (see also ~\cite{F2} for a flow version). An analog statement is true is we switch shadowing by expansive (see \S\ref{stable}  for discrete case and ~\cite{F2} for flows) or by weak specification (see again \S\ref{stable}). Therefore, any volume-preserving diffeomorphism in the $C^1$-interior of the $C^1$-com\-ple\-men\-tary of the volume-preserving Anosov maps must be $C^1$-accumulated by some volume-preserving diffeomorphisms which cannot exhibit the shadowing property, by some volume-preserving non expansive diffeomorphism and also by one without the weak specification property. In other words, $C^1$-far from the volume-preserving Anosov maps we have that systems with the absent of the shadowing property are dense and the same holds for non expansive or  conservative maps without weak specification.  Nevertheless, it was still unknown the ``quantity'' of non shadowable or non expansive maps or non weakly specification maps in the complementary of the Anosov ones and, furthermore, the \emph{non coexistence} of all properties. Our main contribution (cf. Theorems~\ref{thm1},~\ref{thm11} and ~\ref{thm3}) is to answer this question by proving that non shadowable, non expansive, and non weak specification conservative maps are, indeed, a Baire second category or non-meagre subset when considering the $C^1$-complementary of Anosov conservative maps. We also prove, in Theorem~\ref{thm2}, Theorem~\ref{thm22} and Theorem~\ref{thm33}, the corresponding version for flows which are volume-preserving. In conclusion we give a contribution in favor of ~\cite[Conjecture 1]{AD}.

We point out that a manifestation of the phenomenon of genericity of non shadowing away from hyperbolicity was considered in \cite{AD,ALL,LW} for dissipative systems. With respect to the analog results for non expansive ones we refer the results ~\cite{A,Se}. For specification see ~\cite{SSY,ASS}.

We observe that some care is needed to treat the low dimensional case both for flows and maps. That is the content of section \S\ref{SC} and somehow related with Smale conjecture on the $C^1$-density of hyperbolicity on surfaces.

Finally, in \S\ref{stable}, we also prove that the $C^1$-robusteness of the aforesaid properties guarantee that the system is Anosov.

\subsection{Basic definitions for volume-preserving
diffeomorphisms}

Let $M=M^d$ be a $d$-di\-men\-sio\-nal ($d\geq 2$) Riemannian closed and
connected manifold and
 let $d(\cdot,\cdot)$ denotes the distance on $M$ inherited by the Riemannian structure. We endow $M$ with a volume-form (cf. ~\cite{Mo})
  and let $\mu$ denote the Lebesgue measure related to it. Let ${\rm Diff}_\mu^{1}(M)$ denote the set of volume-preserving (or conservative) diffeomorphisms
   defined on $M$, i.e. those diffeomorphisms $f$ such that $\mu(B)=\mu(f(B))$ for any $\mu$-measurable subset $B$.
   Consider this space endowed with the $C^1$ Whitney topology. The Riemannian inner-product induces a norm $\|\cdot\|$
    on the tangent bundle $T_x M$. We will use the usual uniform norm of a bounded linear map $A$ given by $\|A\|=\sup_{\|v\|=1}\|A\cdot v\|$.

Let $f\in {\rm Diff}_\mu^{1}(M)$. Given $\delta>0$, we say that a
sequence of points $\{x_i\}_{i\in\mathbb{Z}}\subset M$ is a
\emph{$\delta$-pseudo-orbit of $f$} if $d(f(x_i),x_{i+1})<\delta$
 for all $i\in\mathbb{Z}$. Let $\Lambda\subseteq M$ be a closed
$f$-invariant set. We say that $f$ has the {\it shadowing
property on $\Lambda$} if for every $\epsilon>0$ there is
$\delta>0$ such that for any $\delta$-pseudo orbit
$\{x_i\}_{i=a}^{b}\subset\Lambda$ of $f$ (where $-\infty\leq a<b\leq
\infty),$ there is a point $y\in M$ such that $d(f^i(y),
x_i)<\epsilon$ for all $a\leq i\leq b-1.$ The diffeomorphism $f$
has the {\it shadowing property} if $\Lambda=M$ in the above
definition. Note that $f$ has the shadowing property if and only
if $f^n$ has the shadowing property for all $n\in\Z.$  

A diffeomorphism $f\colon M\rightarrow M$ is said to be \emph{expansive} (cf. Ma\~n\'e ~\cite{M0}) if for any $x,y\in M$ we have $d(f^n(x),f^n(y))<\epsilon$ for all $n$, then $x=y$. 

The diffeomorphism $f\colon M\rightarrow M$ satisfies the \emph{specification property} (see ~\cite{DGS}) if for all $\epsilon>0$ there exists $N\in\mathbb{Z}^+$, such that for any $k\geq 2$, $ x_1,...x_k\in M$ and $a_1\leq b_1<a_2\leq b_2<...<a_k\leq b_k$, 
with $a_i-b_{i-1}\geq N$ for $i\in\{2,...,k\}$ and $a_i, b_i\in \mathbb{Z}$, there exists $y\in M$ such that:
$$d(f^j(y),f^j(x_i))\leq \epsilon,\,\,\,\,\,\text{for all $j\in\{a_i,...,b_i\}$ and $i\in\{1,...,k\}$}.$$

If we take $k=2$ in the previous definition then we say that the maps satisfies the \emph{weak specification property}.

We say that $\Lambda\subseteq M$ is a {\it transitive set} if there is a point $x
 \in\Lambda$ such that $\omega(x)=\Lambda.$ We say that $\Lambda$ is {\it hyperbolic} if the tangent bundle
$T_{\Lambda}M$ has a $Df$-invariant splitting $E^s\oplus E^u$ and
there exists constants  $C>0$ and $0<\lambda<1$ such that
$$\|Df^n(x)|_{E_x^s}\|\leq C\lambda^n\;\;{\rm and}\;\;\|Df^{-n}(x)|_{E_x^u}\|\leq C\lambda^{n} $$
for all $x\in \Lambda$ and $n\geq 0.$ If $\Lambda=M$, then $f$ is called 
\emph{Anosov}. Let $P_h(f)$ stands for the subset of $M$ defined by the hyperbolic periodic points of $f$ and $P(f)\supseteq P_h(f)$ the subset of all periodic orbits of $f$.

We reached a point where it makes sense to present our main results for the discrete-case. In the following we consider $M=M^d$ with $d\geq 3$.

\begin{maintheorem}\label{thm1}
There exists a $C^1$-residual subset $\mathcal{R}\subset {\rm Diff}_\mu^{1}(M)$ such that $f\in\mathcal{R}$ has the shadowing property if and only
if $f$ is Anosov.
\end{maintheorem}

\begin{maintheorem}\label{thm11}
There exists a $C^1$-residual subset $\mathcal{R}\subset {\rm Diff}_\mu^{1}(M)$ such that $f\in\mathcal{R}$ is expansive if and only
if $f$ is Anosov.
\end{maintheorem}

\begin{maintheorem}\label{thm3}
There exists a $C^1$-residual subset $\mathcal{R}\subset {\rm Diff}_\mu^{1}(M)$ such that $f\in\mathcal{R}$ has the weak specification property if and only if $f$ is Anosov.
\end{maintheorem}

As a direct corollary of previous results we have:

\begin{cor}
There exists a $C^1$-residual subset $\mathcal{R}\subset {\rm Diff}_\mu^{1}(M)$ such that $f\in\mathcal{R}$ is expansive, has the shadowing property and has the weak specification property if and only
if $f$ is Anosov.
\end{cor}

\subsection{Basic definitions for divergence-free vector fields}
When dealing with flows we let $M=M^d$ be a $d$-dimensional ($d\geq 3$) Riemannian closed and
connected manifold.
 Given a $C^r$ ($r\geq 1$) vector field $X\colon M \rightarrow TM$ the solution of the equation $\dot x=X(x)$ gives rise to a $C^r$ flow,
  $X^t$; by the other side given a $C^r$ flow we can define a $C^{r-1}$ vector field by considering $X(x)=\frac{d X^t(x)}{dt}\vert_{t=0}$.
  We say that $X$ is \emph{divergence-free} if its divergence is equal to zero. Note that, by Liouville formula,
   a flow $X^t$ is volume-preserving (or incompressible, or conservative) if and only if the corresponding vector field, $X$, is divergence-free.

Let $\mathfrak{X}_\mu^r(M)$ denote the space of $C^r$
divergence-free vector fields and we consider the usual $C^r$
Whitney topology on this space.

For $\delta>0,$ we say that
$$\{(x_i, t_i):x_i\in M,
t_i\geq1\}_{i\in\Z}$$ is a \emph{$(\delta, 1)$-pseudo orbit} of $X$ if
$d(X^{t_i}(x_i), x_{i+1})<\delta$ for all $i\in\Z$. Define {\sf
Rep} as the set of increasing homeomorphisms $h:\R\to \R$ such
that $h(0)=0.$ Fix $\epsilon>0$ and define {\sf Rep}($\epsilon$)
as follows:
$${\sf Rep}(\epsilon)=\Big\{ h\in{\sf Rep}:
\Big|\frac{h(t)}{t}-1\Big|<\epsilon\Big\}.$$ Let $\Lambda\subseteq
M$ be a compact $X^t$-invariant set. We say that $X^t$ has {\it
the shadowing property} on $\Lambda$ if for any $\epsilon>0$ there
is $\delta>0$ such that for any $(\delta, 1)$-pseudo orbit
$\{(x_i, t_i)\}_{i\in\Z}\subset\Lambda$, let
$T_i=t_0+t_1+\cdots+t_i$ for any $0\leq i<b,$ and
$T_i=-t_{-1}-t_{-2}-\cdots-t_i$ for any $a<i< 0$, there exist a
point $y\in M$ and an increasing homeomorphism $h:\R\to\R$ with
$h(0)=0$ such that $$d(X^{h(t)}(y), X^{t-T_{i}}(x_i))<\epsilon,$$
for any $T_i<t<T_{i+1}.$ If $\Lambda=M$, then $X^t$ has the
shadowing property.

A flow $X^t$ is said to be \emph{expansive} (see~\cite{BW})
if for any $\epsilon>0$ there is $\delta>0$ such that if we have $d(X^t(x),X^{\alpha(t)}(y))\leq\delta$, 
for all $t\in\mathbb{R}$, for $x,y\in M$ and a continuous map $\alpha:\mathbb{R}\rightarrow\mathbb{R}$
with $\alpha(0)=0$, then $y=X^{s}(x)$, where $\left|s\right|\leq\epsilon$. 

A {\em specification} $S=(\tau,P)$ for a flow $X^t$ consists of a finite collection 
$\tau=\{I_1,\cdots,I_m\}$ of bounded intervals $I_i=[a_i,b_i]$
of the real line and a map $P:\cup_{I_i \in \tau}I_i \rightarrow \Lambda$ such that for any $t_1,t_2 \in I_i$ we have
$$X^{t_2}(P(t_1))=X^{t_1}(P(t_2)).$$
A specification $S$ is said to be $K$-{\em spaced} if $a_{i+1} \geq b_i+K$ for all $i \in \{1,\cdots,m\}$ and the minimal 
$K$ is called the spacing of $S$. If $\tau=\{I_1,I_2\}$, then $S$ is said to be a {\em weak specification}.
We say that $S$ is $\epsilon$-{\em shadowed} by $x \in \Lambda$ if 
$d(X^t(x),P(t))<\epsilon$ for all $t \in \cup_{I_i \in \tau}I_i$.

We say that $M$ has the {\em weak specification property} if for any $\epsilon>0$ 
there exists a $K=K(\epsilon) \in \Rr$ such that any $K$-spaced weak specification $S$ is
$\epsilon$-shadowed by a point of $M$.

Given a vector field $X$ we denote by $\text{Sing}(X)$ the set of
\emph{singularities} of $X$, i.e. those points $x\in M$ such that
$X(x)=\vec{0}$. Let $M':=M\setminus \text{Sing}(X)$ be the set of
\emph{regular} points. Given $x\in M'$ we consider its normal
bundle $N_{x}=X(x)^{\perp}\subset T_{x}M$ and define the
associated \emph{linear Poincar\'{e} flow} by
$P_{X}^{t}(x):=\Pi_{X^{t}(x)}\circ DX^{t}({x})$ where
$\Pi_{X^{t}(x)}\colon T_{X^{t}(x)}M\rightarrow N_{X^{t}(x)}$ is
the projection along the direction of $X(X^{t}(x))$. In the same
way as we did in the discrete-time case we can define
uniform hyperbolicity, dominated splitting, partial hyperbolicity
and volume-hyperbolicity for the linear Poincar\'e flow $P_X^{t}$,
in subsets of $M'$ and related to subbundles of the normal bundle
$N$. We also observe that, when $\Lambda\subseteq M$ is compact,
the hyperbolicity of the tangent map $DX^t$ on $\Lambda$
implies the hyperbolicity for the linear Poincar\'e flow
of $X$ on $\Lambda$ (see \cite[Proposition 1.1]{D}).

At this moment we present the corresponding versions of Theorems ~\ref{thm1}, ~\ref{thm11} and ~\ref{thm3} for flows considering $M=M^d$ with $d\geq 4$.

\begin{maintheorem}\label{thm2}
There exists a $C^1$-residual subset $\mathcal{R}\subset \mathfrak{X}_\mu^r(M)$ such that $X\in\mathcal{R}$ has the shadowing property if and only
if $X$ is Anosov\footnote{We would like to thank Raquel Ribeiro for pointing us out that she also proved independently, in~\cite{Rib}, this same result.}.
\end{maintheorem}

\begin{maintheorem}\label{thm22}
There exists a $C^1$-residual subset $\mathcal{R}\subset \mathfrak{X}_\mu^r(M)$ such that $X\in\mathcal{R}$ has the expansive property if and only
if $X$ is Anosov.
\end{maintheorem}

\begin{maintheorem}\label{thm33}
There exists a $C^1$-residual subset $\mathcal{R}\subset {\rm Diff}_\mu^{1}(M)$ such that $f\in\mathcal{R}$ has the weak specification if and only
if $X$ is Anosov.
\end{maintheorem}

The proof of Theorem~\ref{thm33}  has been sketched in ~\cite[Theorem 9.1]{ASS}. Here we give a slightly different and more direct proof.

Again, as a trivial conclusion we have:
\begin{cor}
There exists a $C^1$-residual subset $\mathcal{R}\subset \mathfrak{X}_\mu^r(M)$ such that $X\in\mathcal{R}$ is expansive has the shadowing  property and has the weak specification property if and only if $X$ is Anosov.

\end{cor}

\section{Proof of Theorems \ref{thm1}, \ref{thm11} and \ref{thm3}}

\subsection{Homoclinic points, index and weak hyperbolicity}

The strategy that we are going to use to obtain our results is a \emph{change of index}  argument on hyperbolic periodic points and its link to non hyperbolic behavior.

It is well known from invariant manifolds theory that, if $p\in P_h(f)$ and has period $\pi$, then the sets
\begin{align*}W^s(p,f)=\{x\in M: f^{\pi n}(x)\to p\quad\mbox{as}\quad
n\to\infty\}\quad\mbox{and}\\W^u(p,f)=\{x\in M: f^{-\pi n}(x)\to
p\quad\mbox{as}\quad n\to\infty\}\end{align*} are
$C^1$-injectively immersed submanifolds of $M.$ Let $q\in P_h(f)$, we say that $p$ and $q$ are {\it
homoclinically related}, and write $p\sim q$ if $$W^s(p)\pitchfork
W^u(q)\not=\emptyset\quad\mbox{and}\quad W^u(p)\pitchfork
W^s(q)\not=\emptyset,$$
where $\pitchfork$ stands for the transversal intersection. It is clear that if $p\sim q$, then ${\rm
index}(p)={\rm index}(q)$; i.e., ${\rm dim}W^s(p)={\rm
dim}W^s(q)$. When there is no ambiguity we write $W^\sigma(p)$ instead of $W^\sigma(p,f)$ for $\sigma=u,s$.

In \cite{BC} Bonatti and Crovisier proved the following quite useful and core result on $C^1$-generic conservative dynamics.

\begin{thm}(\cite[Th\'eor\`eme 1.3]{BC})\label{transitive} There is a residual
set $\mathcal{R}_1\subset{\rm Diff}_{\mu}^1(M)$ such that for any
$f\in\mathcal{R}_1$, $f$ is transitive.
\end{thm}

Note that, in the volume-preserving context, being a map transitive or not, it does not have sinks nor
sources. Moreover, hyperbolic periodic orbits must be of saddle-type.

For a given $\delta>0$, we say that $p\in
P(f)$ has a {\it $\delta$-weak eigenvalue} if there is an
eigenvalue $\lambda$ of $Df^{\pi(p)}(p)$ such that
$(1-\delta)^{\pi(p)}<|\lambda|<(1+\delta)^{\pi(p)},$ where
$\pi(p)$ is the period of $p$.

\subsection{Shadowing} The following result allows us to obtain, in a $C^1$-generic sense, connections between invariant manifolds from the shadowing property.

\begin{lem}\label{tan} Let $f\in\mathcal{R}_1,$ and let $p, q\in P_h(f)$. If $f$ has the shadowing property, then $W^s(p)\cap
W^u(q)\not=\emptyset$.
\end{lem}
\noindent{\bf Proof.} By conservativeness, we have that $p, q\in P_h(f)$ are saddles. Since
$f\in\mathcal{R}_1,$ $f$ is transitive. Then, we can find a point
$x\in M$ such that $\omega_f(x)=M.$ Since $p, q$ are hyperbolic
saddles, we take  $\epsilon(p)>0$ and $\epsilon(q)>0$  as before
with respect to $p$ and $q.$  Fix $\epsilon={\rm
min}\{\epsilon(p), \epsilon(q)\}.$ To simplify the notation in the
proof, assume that $f(p)=p$ and $f(q)=q.$ Let
$0<\delta<\epsilon/2$ be the number of the shadowing property of
$f.$ For the chosen $\delta>0$, we can find $\ell_1>0,$ and $\ell_2>0$ such
that $d(f^{\ell_1}(x), p)<\delta$ and $d(f^{\ell_2}(x), q)<\delta.$ We
may assume that $\ell_2>\ell_1.$ Then, we can define a finite
$\delta$-pseudo orbit $\{p, f^{\ell_1}(x), f^{\ell_1+1}(x), \ldots,
f^{\ell_1+\ell_2-1}(x), q\}$. Let us now extend the finite $\delta$-pseudo
orbit as follows: Put $x_i=f^i(p)$ for all $i\leq \ell_1,$
$x_i=f^{\ell_1+i}(x)$ for all $0\leq i\leq\ell_2-1,$ and $x_i=f^i(q)$
for all $i\geq \ell_2.$ Also, we set, via a time translation, the $\delta$-pseudo orbit again
as follows: let $\ell_1=0$ and $\ell_2=\ell.$ Then $x_i=f^i(p)$ for
$i\leq 0,$ $x_i=f^i(x)$ for $0<i<\ell-1$ and $x_{\ell+i}=f^{i}(q)$ for
$i\geq \ell.$

Therefore, the $\delta$-pseudo orbit is
\begin{align*}\{x_i\}_{i\in\Z}&=\{\ldots, p, f^{\ell_1}(x),
f^{\ell_1+1}(x), \ldots, f^{\ell_1+\ell_2-1}(x), q, \ldots\}\\&=\{\ldots,
x_{-1}, x_0(=p), x_1, x_2, \ldots, x_{\ell-1}, x_\ell(=q), x_{\ell+1},
\ldots\}.
\end{align*}

By the shadowing property, there is a point $y\in M$ such that
$d(f^i(y), x_i)<\epsilon$ for all $\in\Z.$ Then for all $i\leq0$
$$d(f^i(y), x_i)=d(f^i(y), p)<\epsilon,$$ and for all $i\geq\ell$,
$$d(f^i(y), x_i)=d(f^i(y), q)<\epsilon.$$ 

In conclusion, we see that $y\in
W^u_{\epsilon}(p)$ and $f^\ell(y)\in W^s_{\epsilon}(q),$ and so,
$y\in W^u(p)\cap W^s(q).$ The other case is similar. Thus, we obtain
$W^s(p)\cap W^u(q)\not=\emptyset$ and $W^u(p)\cap
W^s(q)\not=\emptyset$ and the lemma is proved. $\square$\\

The next result shows that $C^1$-generically shadowable maps are \emph{homogeneous}, i.e., their periodic orbits have constant index.

\begin{lem}\label{transe} There is a residual set
$\mathcal{R}_2\subset{\rm Diff}_{\mu}^1(M)$ such that for any
$f\in\mathcal{R}_2$, if $f$ has the shadowing property, 
 then for any $p,q\in P_h(f)$ we have ${\rm index}(p)={\rm index}(q)$.
\end{lem}

\noindent {\bf Proof.} The result is meaningless if $\dim M=2$ (see \S\ref{SWS}). Assume that $\dim M\geq 3$. Let $\mathcal{R}_2=\mathcal{R}_1\cap \mathcal{K}$, where $\mathcal{K}$ are the Kupka-Smale volume-preserving diffeomorphisms given by Robinson's theorem (see~\cite{R}), and let $f\in\mathcal{R}_2$. Suppose that $f$ has the
shadowing property. By the conservative constraint $p, q\in P_h(f)$ are saddles. The proof is by contradiction, so let us assume that ${\rm index}(p)\not={\rm
index}(q).$ Then ${\rm dim}W^s(p)+{\rm dim}W^u(q)<{\rm dim}M$ or
${\rm dim}W^u(p)+{\rm dim}W^s(q)<{\rm dim}M.$ Without loss of
generality, we may assume that ${\rm dim}W^s(p)+{\rm
dim}W^u(q)<{\rm dim}M$. Since $f$ is Kupka-Smale we have $W^s(p)\cap
W^u(q)=\emptyset$ which, by Lemma \ref{tan}, is a contradiction.
$\square$\\

\medskip

The following proposition is the discrete and conservative version of \cite[Lemma 5.1]{ASS} and its proof is completely analog. 

\begin{pro} \label{ASS2}
There exists a residual subset $\mathcal{R}_3\subset{\rm Diff}_{\mu}^1(M)$ such that if $f\in\mathcal{R}_3$ is $C^1$-approximated by $\{f_n\}_{n\in\mathbb{N}}$ such
that each $f_n\in{\rm Diff}_{\mu}^1(M)$ has two distinct points, $x_n,y_n\in P_h(f_n)$, with different indices and
with $d(x_n, y_n) < \epsilon$, then there exist two distinct points, $x, y \in P_h(f)$, with different indices
and with $d(x, y) < 2\epsilon$.
\end{pro}

The next result shows how to create, via $C^1$-perturbations and within the conservative setting, hyperbolic periodic points close to each other and with different indices once we have a non hyperbolic periodic orbit.

\begin{lem}\label{PF}
Let $f\in {\rm Diff}_{\mu}^1(M^d)$ with $d\geq 3$. If $p\in P(f)\setminus P_h(f)$ and $\epsilon>0$, then there is $g\in {\rm Diff}_{\mu}^1(M)$ arbitrarily 
$C^1$-close to $f$ such that $g$ has two periodic
orbits $p_1, p_2\in P_h(g)$ with different indices and with $d(p_1,p_2)<\epsilon$.
\end{lem}

\noindent {\bf Proof.} 
Let $p\in P(f)$ be the non-hyperbolic periodic orbit of period $\pi$ and $\epsilon>0$. We use Pugh-Robinson's closing lemma for volume-preserving diffeomorphisms (\cite[\S 8 (b)]{PR}) to obtain $f_1\in\text{Diff}^1_\mu(M)$, such that $f_1$ is  arbitrarily $C^1$-close to $f$, with $p_2\in P(f_1)$ and $\epsilon$-close to $p$ by closing some recurrent orbit\footnote{By Poincar\'e recurrence almost every point is recurrent.}. Moreover, since hyperbolicity holds open and densely even in the volume-preserving setting (recall that $d\geq 3$), $p_2$  can be chosen hyperbolic. Denote $\text{index}(p_2)=i$. After this perturbation  we get $f_1\in \text{Diff}^1_\mu(M)$ such that $f_1$ has a periodic orbit $p_1$, say $\frac{\epsilon}{10}$-close to $p$, with period $\pi$. We observe that $p_1$ may not be the analytic continuation of $p$ and this is precisely the case when $1$ is an eigenvalue of the tangent map $Df^{\pi}(p)$. If $p_1$ is not hyperbolic take $f_2=f_1$. If $p_1$ is hyperbolic for $Df_1^{\pi}(p_1)$, then, since $f_1$ is arbitrarily $C^1$-close to $f$, the distance between the spectrum of $Df_1^{\pi}(p_1)$ and the unitary circle can be taken arbitrarily close to zero. This means that we are in the presence of a very weak hyperbolicity, i.e., of a $\delta$-weak eigenvalue (for some small $\delta>0$), thus in a position to apply ~\cite[Proposition 7.4]{BDP} to obtain $f_2\in \text{Diff}^1_\mu(M)$, such that $p_1$ is  a non-hyperbolic periodic orbit. Moreover, this \emph{local} perturbation can be done far from the hyperbolic periodic point $p_2$. Once again we use \cite[Proposition 7.4]{BDP} in order to obtain $g\in \text{Diff}^1_\mu(M)$ such that $p_1\in P(g)$, is hyperbolic and $\text{index}(p_1)\not= i$ and still $d(p_1,p_2)<\epsilon$.
 $\square$\\

\begin{lem}\label{weak} There is a residual set
$\mathcal{R}_4\subset{\rm Diff}_{\mu}^1(M)$ such that for any
$f\in\mathcal{R}_4$, if $f$ has the shadowing property  then there
is $\delta>0$ such that for any $p\in P(f),$ $p$ does not have a
$\delta$-weak eigenvalue.
\end{lem}

\noindent {\bf Proof.} Let $\mathcal{R}_4=\mathcal{R}_2\cap\mathcal{R}_3,$
 and let assume that $f\in\mathcal{R}_4$ satisfies the shadowing property. The proof is by contradiction. Suppose
 that there is $p\in P(f)$ such that for any $\delta>0$,
 $p$ has a $\delta$-weak eigenvalue.  Thus, by Lemma~\ref{PF}, we can
 find $g\in {\rm Diff}_{\mu}^1(M)$, $C^1$-close to $f$, such that $g$ has two hyperbolic
 periodic orbits $p_1, p_2\in P(g)$ with ${\rm
 index}(p_1)\not={\rm index}(p_2)$. But $f\in \mathcal{R}_2$ and this contradicts Lemma \ref{transe}.$\square$\\

\begin{rk}\label{rmkexpansive}
Previous result also holds if we replace \emph{shadowing property} by the \emph{expansive property}. Actually, the proof in ~\cite[Lemma 5.3]{A}  can be followed step by step and we just have to use Franks' lemma for conservative maps proved in ~\cite[Proposition 7.4]{BDP}.
\end{rk}

With respect to $\delta$-weak eigenvalues, we have the following lemma which was proved by Arbieto and works equally within conservative constraints.
\begin{lem}\label{two}{\em\cite[Lemma 5.1 (2)]{A}} There is a residual set
$\mathcal{R}_5\subset {\rm Diff}_{\mu}^1(M)$ such that for any
$f\in\mathcal{R}_5$ the following holds:  fixing $\delta>0$, if for any $C^1$-neighborhood
$\mathcal{U}(f)$, there is $g\in\mathcal{U}(f)$ and 
$p_g\in P_h(g)$, such that $p_g$ has a $\delta$-weak eigenvalue, then there exists $p\in P_h(f)$
with a $2\delta$-weak eigenvalue.
\end{lem}

\noindent
 {\bf Proof of Theorem ~\ref{thm1}.} 
 
  We define the set $\mathcal{F}_{\mu}^\star(M)$ as
the set of \emph{star volume-preserving diffeomorphisms}, i.e., those diffeomorphisms $f\in{\rm Diff}_{\mu}^1(M)$ which have a
$C^1$-neighborhood $\mathcal{U}(f)\subset{\rm Diff}_{\mu}^1(M)$ such
that for any $g\in\mathcal{U}(f)$, every periodic point of $g$ is
hyperbolic. 

 We begin to prove that there is a residual set
$\mathcal{R}'\subset{\rm Diff}_{\mu}^1(M)$ such that for any
$f\in\mathcal{R}',$ if $f$ has the shadowing property, then
$f\in\mathcal{F}_{\mu}^\star(M)$, and so $f$ is Anosov (Theorem~\ref{Ar}).
 
 Let
 $\mathcal{R}'=\mathcal{R}_4\cap\mathcal{R}_5$, and let
 $f\in\mathcal{R}'$ satisfying the shadowing property. To derive a contradiction, we may assume that
 $f\not\in\mathcal{F}_{\mu}^\star(M)$. Then, for any arbitrarily small $C^1$-neighborhood $\mathcal{U}(f)\subset{\rm Diff}_{\mu}^1(M)$ there exists $g\in\mathcal{U}(f)$, with a non-hyperbolic periodic orbit. Hence, there is a periodic point $p_g\in P(g)$ such
 that $p_g$ has a
 $\delta/2$-weak eigenvalue. By Lemma \ref{two}, there exists $p\in P(f)$ with a $\delta$-weak
 eigenvalue. But by Lemma \ref{weak} we reach a contradiction. Thus we must have $f$ in $\mathcal{F}_{\mu}^\star(M)$ and Lemma~\ref{Ar} concludes the proof. $\square$\\

\subsection{Expansiveness}

In this section we present the proof of Theorem~\ref{thm11} taking into account several results previously obtained.

\medskip

 {\bf Proof of Theorem \ref{thm11}.} We define the residual subset of the theorem by the intersection of $\mathcal{R}_5$ with the residual given by Remark~\ref{rmkexpansive}. The proof is obtained by \emph{reductio ad absurdum}. Let $f$ be expansive and in the previous defined residual. To derive a contradiction, we again assume that
 $f\not\in\mathcal{F}_{\mu}^\star(M)$ and proceed as in the proof of Theorem~\ref{thm1}. We leave the details to the reader. $\square$\\

\subsection{Weak specification} 

The proof of the following result is the same as the one in \cite[Lemma 2.1]{SSY} since this lemma is abstract and holds also in the volume-preserving case.

\begin{lem}\label{manifolds}
Let $f\in \text{Diff}^{\,\,1}_\mu(M)$. If $q_1,q_2\in P_h(f)$ and $f$ has the weak specification property, then $W^u(q_1,f)\cap W^s(q_2,f)\not=\emptyset$.
\end{lem}

The following result will be the essential key to obtain the proof of Theorem~\ref{thm3}.

\begin{lem}\label{spec}
There is a residual set
$\mathcal{R}_6\subset{\rm Diff}_{\mu}^1(M)$ such that, for any
$f\in\mathcal{R}_6$, if $f$ has the weak shadowing property, then there
is $\delta>0$ such that for any $p\in P(f),$ $p$ does not have a
$\delta$-weak eigenvalue.
\end{lem}
 {\bf Proof.} We define $\mathcal{R}_6=\mathcal{K}\cap\mathcal{R}_3$, where $\mathcal{R}_3$ is the $C^1$-residual given in Proposition ~\ref{ASS2} and $\mathcal{K}$ is the residual given by Robinson's Kupka-Smale theorem~\cite{R}. Assume, by contradiction, that there exists $p_n\in P(f)$, such that $p_n$ has a $\frac{1}{n}$-weak eigenvalue. Then, there exists $f_n\in{\rm Diff}_{\mu}^1(M)$ such that $f_n\rightarrow f$, in the $C^1$-sense, displaying $q_1^n,q_2^n\in P(f_n)$ with different indices and such that $d(q_1^n,q_2^n)<\epsilon$. Then, using Proposition~\ref{ASS2}, there exist two distinct hyperbolic periodic points, $q_1, q_2 \in P_h(f)$, with different indices and with $d(x, y) < 2\epsilon$. By Lemma ~\ref{manifolds} $W^u(q_1,f)\cap W^s(q_2,f)\not=\emptyset$ and this intersection cannot be transversal which contradicts the fact that $f\in \mathcal{K}$.
$\square$\\

\noindent
 {\bf Proof of Theorem \ref{thm3}.} We define the residual subset of the theorem by the intersection of $\mathcal{R}_5$ with the residual given by Lemma~\ref{spec}. Once again the proof is obtained by \emph{reductio ad absurdum}. Let $f$ be a map with the weak specification property  and in the previous defined residual. To derive a contradiction, we may assume that
 $f\not\in\mathcal{F}_{\mu}^\star(M).$ Then there is a periodic point $p$ such
 that $p$ has a
 $\delta/2$-weak eigenvalue.  By Lemma \ref{two}, $p$ has a $\delta$-weak
 eigenvalue which contradicts Lemma~\ref{spec} and Theorem~\ref{thm3} is proved. $\square$\\

\section{The low dimensional case}\label{SC}

\subsection{Shadowing and weak specification}\label{SWS} With respect to the shadowing and weak specification properties we were not able to obtain the corresponding version of our main theorems. Next, we discuss shortly the main difficulties in order to reach the same conclusions as in the higher dimensional case.

We observe that the change of index argument makes no sense within the area-preserving case; the index cannot be two, by conservative reasoning, hence it is always equal to one. Furthermore, the correspondent to the non hyperbolic phenomena related to the change of index in lower dimension is the mechanism of homoclinic tangencies.

The natural approach now is to obtain an adapted version of Proposition ~\ref{ASS2}, namely, prove that there exists a residual subset $\mathcal{R}\subset{\rm Diff}_{\mu}^1(M^d)$ ($d=2$) such that if $f\in\mathcal{R}$ is $C^1$-approximated by $\{f_n\}_{n\in\mathbb{N}}$ such
that each $f_n\in{\rm Diff}_{\mu}^1(M)$ has $x_n\in P_h(f_n)$, with a homoclinic tangency $\tau_n$ and
with $d(x_n, \tau_n) < \epsilon$, then there exists $x \in P_h(f)$, with a homoclinic tangency $\tau$
and with $d(x, \tau) < 2\epsilon$. However, in order to go on with the generic argument we strongly need the $C^1$-\emph{robusteness} of homoclinic tangencies which is still unknown and directly related to Smale conjecture on the $C^1$-density of hyperbolicity on surfaces (see the discussion on ~\cite[\S 1.3]{BoD}). Finally, taking into account recent results in ~\cite{Mor}, we conjecture that:

\begin{con}
For non Anosov area-preserving maps the non-shadowable (or non weak specification) ones are meagre in the $C^1$-sense.
\end{con}

We believe that the same result holds for incompressible flows on three-manifolds.

\subsection{Expansiveness}  
Fortunately, when dealing with the expansiveness property we were able to obtain a proof. We begin by presenting some useful definitions.
We recall that a periodic point is called \emph{elementary} if, for the eigenvalues $\sigma_i$ (on the period),  the following holds; if $\prod_{p_i\in\mathbb{Z}}\sigma_i^{p_i}=1$, then $p_i=0$ for all $i$. That is, its eigenvalues are multiplicatively independent over $\mathbb{Z}$. Elementary points have simple spectrum and none of its eigenvalues is a root of unity or equals to 1. In our conservative maps case, if $\dim M=2$, then we have three possibilities: a \emph{hyperbolic saddle} (real eigenvalues $\sigma_2=\sigma_1^{-1}$), an \emph{elliptic point} (non real eigenvalues, conjugated and of norm one) and a \emph{parabolic point} (with eigenvalues equal to 1 or $-1$). Observe that the first two cases are robust under small perturbations. Moreover, elementary elliptic points are associated to an irrational rotation number. By ~\cite[Theorem 1 B. i)]{R}, if  $\dim M=2$, there exists a residual subset in ${\rm Diff}_{\mu}^1(M)$ such that any element in this residual displays all its elliptic points of elementary type.

In order to obtain the proof of Theorem~\ref{thm11} on surfaces we begin by proving the following preliminary result. 

\begin{lem}\label{first} If $f\in{\rm Diff}_{\mu}^1(M)$ has a non-hyperbolic point $p$, $\dim M=2$ and $\epsilon,\nu>0$, then there exists $g\in{\rm Diff}_{\mu}^1(M)$, such that:
\begin{enumerate}
\item $g$ is $\epsilon$-$C^1$-close to $f$;
\item $g$ has two hyperbolic periodic points $q_1$ and $q_2$ with the same period and
\item $\underset{i\in\mathbb{Z}}{\sup}\,d(g^i(q_1),g^i(q_2))<\nu$.
\end{enumerate}
\end{lem}

\noindent {\bf Proof.} Let $p$ be a non-hyperbolic period orbit for $f$ with period $\pi$. Non-hyperbolic points are elliptic or parabolic. Since ellipticity is stable and any parabolic point can be made elliptic by a small perturbation we choose $g_1\in{\rm Diff}_{\mu}^\infty(M)$ $\epsilon/4$-$C^1$-close to $f$ given by Zehnder's smoothing theorem (see~\cite{Z}) such that $p$ is elliptic w.r.t. $g_1$. Now, using Franks' lemma (\cite[Proposition 7.4]{BDP}) we consider $g_2\in{\rm Diff}_{\mu}^\infty(M)$ $\epsilon/4$-$C^1$-close to $g_1$ such that, $Dg_2^{\pi}(p)$ is a rational rotation. Then, there exists a positive integer $n$ such that $Dg_2^{n\pi}(p)=id$. Now applying the Pasting Lemma (\cite{AM}) along the orbit of $p$, we
can get an $\epsilon/4$-$C^1$-perturbation $g_3\in{\rm Diff}_{\mu}^\infty(M)$ of $g_2$ such that there is a small
ball $B(p,\nu/4)$ satisfying the equality $g_3^{n\pi}|_{B(p,\nu/4)}=id$. Hence, we get two periodic
point $q_1$ and $q_2$ in $B(p,\nu/4)$ with the same period such that ${\sup}_{i\in\mathbb{Z}}\,d(g^i(q_1),g^i(q_2))<\nu/2$. By making an $\epsilon/4$-$C^1$-perturbation $g\in{\rm Diff}_{\mu}^\infty(M)$ we can turn $q_1$ and $q_2$ hyperbolic and such that ${\sup}_{i\in\mathbb{Z}}\,d(g^i(q_1),g^i(q_2))<\nu$. This ends the proof of the lemma.
 $\square$\\

Note that by a result by Newhouse (\cite{N}), if
$f\in\mathcal{F}_{\mu}^\star(M)$ and $\dim M=2$, then $f$ is Anosov. Recently, Arbieto and Catalan \cite{AC} proved that, in any dimension, 
 $f\in\mathcal{F}_{\mu}^\star(M)$ is
Anosov (Theorem~\ref{Ar}). To prove this, they used the Ma\~n\'e's results in
\cite[Proposition II.1]{M}  and showed that $\overline{P(f)}$ is
hyperbolic if $f\in\mathcal{F}_{\mu}^\star(M)$.
 
\begin{thm}\label{Ar} {\em \cite[Theorem 1.1]{AC}} Every diffeomorphism in
$\mathcal{F}_{\mu}^\star(M)$ is Anosov.
\end{thm} 
 
Since we are in the setting of surfaces the next useful result only needs Newhouse's two-dimensional version of previous result (\cite{N}).
 
\begin{lem}\label{second} There is a residual set
$\mathcal{R}\subset{\rm Diff}_{\mu}^1(M)$, where $\dim M=2$, such that for any
$f\in\mathcal{R}$ outside $\mathcal{F}^\star_\mu(M)$, $f$ has two sequences of periodic points $\{p_n\}_n$ and $\{q_n\}_n$ displaying the same period for each $n$ and
$$\underset{n\rightarrow+\infty}{\lim} \underset{i\in\mathbb{Z}}{\sup} \,d(f^i(p_n),f^i(q_n))=0.$$
\end{lem}

\noindent {\bf Proof.} 
For each positive integer $n$, we denote by $\mathcal{N}_n$ the subset of ${\rm Diff}_{\mu}^1(M)$ such that any $f\in\mathcal{N}_n$  has a $C^1$-neighborhood $\mathcal{U}$ in
${\rm Diff}_{\mu}^1(M)$ with the following $C^1$-open property: for every  $g\in\mathcal{U}$, there are hyperbolic periodic points $p$, $q$ of $g$ with the same period such that 
$$\underset{i\in\mathbb{Z}}{\sup}\,d(g^i(p),g^i(q))<\frac{1}{n}.$$
Let $\mathcal{H}_n$ be the $C^1$-complementary of the $C^1$-closure of $\mathcal{N}_n$. Clearly, $\mathcal{N}_n\cup \mathcal{H}_n$ is $C^1$-open and $C^1$-dense in ${\rm Diff}_{\mu}^1(M)$. We define the $C^1$-residual subset of the lemma by 
$$\mathcal{R}=\underset{n\in\mathbb{N}}{\bigcap}(\mathcal{H}_n\cup \mathcal{N}_n).$$
Let $f\in \mathcal{R}$. If $f\notin \mathcal{F}^\star_\mu(M)$ then
there is a sequence of diffeomorphisms $f_k$ converging to $f$ and a sequence of
non-hyperbolic periodic orbits $\tilde{p}_k$ of $f_k$. Then, for any positive integer $n$, by
the Lemma ~\ref{first}, we have $f\notin \mathcal{H}_n$, and so $f \in \mathcal{N}_n$. 
, for each $n$, $f$ has a $C^1$-neighborhood $\mathcal{U}$ in
${\rm Diff}_{\mu}^1(M)$ with the following property: for all $g\in\mathcal{U}$ (in particular for $f$), there are hyperbolic periodic points $p_n$, $q_n$ for $g$ such that 
$$\underset{i\in\mathbb{Z}}{\sup}\,d(g^i(p_n),g^i(q_n))<\frac{1}{n}.$$

\noindent In conclusion, we can define two sequences of periodic orbits for $f$, $\{p_n\}_n$  and $\{q_n\}$ such that
$$\underset{i\in\mathbb{Z}}{\sup}\,d(f^i(p_n),f^i(q_n))<\frac{1}{n},$$
and so
$$\underset{n\rightarrow+\infty}{\lim} \underset{i\in\mathbb{Z}}{\sup} \,d(f^i(p_n),f^i(q_n))=0.$$
And the proof of the lemma is complete.
 $\square$\\

From the lemma above, we know that if $f \in \mathcal{R}$ and $f$ has the expansive
property, $f$ must be in $\mathcal{F}_{\mu}^\star(M)$. From Theorem~\ref{Ar},
we know that $f$ is Anosov. This ends the proof of Theorem ~\ref{thm11} when $\dim M=2$.

\bigskip

Analogously, we can also obtain the proof of Theorem \ref{thm22} for three-dimensional incompressible flows. The main difference, which is not substantially difficult, is that we have to consider also singularities on the adaptation of Lemma~\ref{first}. Moreover, we are able to use the available versions of the Zehnder's smoothing result (see ~\cite[Theorem 2.2]{AM}), the $C^{1+}$-pasting lemma for vector fields ~\cite[Theorem 3.1]{AM} and also the Franks lemma (~\cite[Lemma 3.2]{BR0}). Finally, in order to obtain a version of Lemma~\ref{second} we make use of the result in  \cite{BR1}.

\section{Proofs of Theorems \ref{thm2}, \ref{thm22} and \ref{thm33}}

\subsection{The case when $d\geq 4$} 

The following result was proved by the first author in \cite{B}.
\begin{thm} There is a residual set $\mathscr{R}_1\subset
\mathfrak{X}_{\mu}^1(M)$ such that if $X\in\mathscr{R}_1,$ then $X$
is a topologically mixing vector field.
\end{thm}

Next, we obtain the continuous-time version of Lemma ~\ref{tan}.

\begin{lem}\label{tan1} Let $X\in\mathscr{R}_1$, and let $\gamma, \eta\in P_h(X)$ be hyperbolic periodic
orbits. If $X$ has the shadowing property, then $W^s(\gamma)\cap
W^u(\eta)\not=\emptyset$.
\end{lem}
\noindent {\bf Proof.} Let $X\in\mathscr{R}_1$, and let $\gamma\in
P_h(X)$ and $\eta\in P_h(X).$ Since $X$ is topologically mixing it is transitive and so there exists a
point $x\in M$ such that $\omega_X(x)=M.$ Since $\gamma, \eta$ are
hyperbolic saddles, we can take $\epsilon(\gamma)>0$ and
$\epsilon(\eta)>0$ such that
$W^{\sigma}_{\epsilon(\gamma)}(\gamma)$, 
$W^{\sigma}_{\epsilon(\eta)}(\eta)$ are as before with respect to
$\gamma$ and $\eta,$ where $\sigma=s, u.$ Take
$\epsilon=\min\{\epsilon(\gamma), \epsilon(\eta)\}.$ Let
$0<\delta=\delta(\epsilon)<\epsilon$ be the number of the
shadowing property of $X.$ Since $X$ is transitive, there are
$k>0$ and $l>0$ such that $d(X^{t_k}(x), \gamma)<\delta$ and
$d(X^{t_l}(x), \eta)<\delta.$ Then we can find $p\in\gamma$ and
$q\in \eta$ such that $d(X^{t_k}(x), p)<\delta$ and $d(X^{t_l}(x),
q)<\delta.$ We may assume that $l>k$. Then we can construct a
finite $(\delta, 1)$-pseudo orbit $\{(x_i, t_i)\}_{i=0}^l$ as
follows:
\begin{itemize}
\item[(1)] $p=x_0,$
\item[(2)] $x_i=X^{t_{k+i}}(x)$ for $t_i\geq 1,$ and $0<i<l-1,$ and
\item[(3)] $x_l=q.$
\end{itemize}
Then $\{(x_i, t_i)\}_{i=0}^l=\{p(=x_0), x_1, x_2, \ldots, x_{l-1},x_l(=q)\}.$ Now, we extend the finite $(\delta, 1)$-pseudo orbit
as follows:
\begin{itemize}
\item[(1)]$X^i(p)=x_i$ for $t_i=1,$ and $i\leq0,$
\item[(2)]$x_i=X^{t_{k+i}}(x)$ for $t_i\geq 1,$ and $0<i<l-1,$ and
\item[(3)]$X^i(q)=x_{l+i}$ for $t_i=1$ and $i\geq0.$
\end{itemize}
Therefore, we have the $(\delta, 1)$-pseudo orbit $$\{(x_i,
t_i)\}_{i\in\Z}=\{\ldots, x_{-1}, x_0(=p), x_1, \ldots, x_{l-1},
x_l(=q), x_{l+1}, \ldots \}.$$
 By the shadowing property, there exist a point $y\in M$ and an
 increasing homeomorphism $h:\R\to \R$ with $h(0)=0$ such that
 $d(X^{h(t)}(y), X^{t-T_i}(x_i))<\epsilon$ for all $i\in\Z,$ and $T_i<t<T_{i+1}.$ Then for
 $i\leq0,$ $X^i(y)\in W^u_{\epsilon}(\gamma)$ and for $t>T_l,$ $X^{t_l}(y)\in
 W^s_{\epsilon}(\eta)$. Thus, we know that $y\in W^u(\gamma)\cap
 W^s(\eta),$ and so $W^u(\gamma)\cap W^s(\eta)\not=\emptyset$. $\square$\\

A vector field $X\in\mathfrak{X}_{\mu}^1(M)$ is said to be {\it
Kupka-Smale} if any element of $Crit(X)$ is hyperbolic and its
invariant manifolds intersect transversely. In \cite{R}, Robinson proved that the set
of Kupka-Smale divergence-vector fields is a $C^1$-residual subset
of $\mathfrak{X}^1_{\mu}(M).$ We denote $\mathscr{K}$ by
the Kupka-Smale $C^1$-residual set in $\mathfrak{X}_{\mu}^1(M)$.

\medskip

The following result establishes the flows version of Lemma~\ref{transe}.

\begin{lem}\label{same} There is a residual set $\mathscr{R}_2\subset
\mathfrak{X}_{\mu}^1(M)$ such that for any $X\in\mathscr{R}_2$, if $X$
has the shadowing property, then all its hyperbolic periodic
orbits have the same index.
\end{lem}
\noindent {\bf Proof.} Let $X\in\mathscr{R}_2=\mathscr{R}_1\cap
\mathscr{K}$ and let $X$ satisfies the shadowing
property. The proof is by contradiction. Let us assume that for some $\gamma,
\eta\in P_h(X)$  we have ${\rm index}(\gamma)\not={\rm index}(\eta).$ Therefore, 
we know that ${\rm dim}W^s(\gamma)+{\rm dim}W^u(\eta)\leq {\rm
dim}M$ or ${\rm dim}W^u(\gamma)+{\rm dim}W^s(\eta)\leq{\rm dim}M.$
Assume that ${\rm dim}W^s(\gamma)+{\rm dim}W^u(\eta)\leq{\rm
dim}M.$ Then we have two cases, (i) ${\rm dim}W^s(\gamma)+{\rm
dim}W^u(\eta)<{\rm dim}M$ or (ii) ${\rm dim}W^s(\gamma)+{\rm
dim}W^u(\eta)={\rm dim}M.$

In case (i) and since $X\in\mathscr{R}_2$, we have that 
$W^s(\gamma)\cap W^u(\eta)=\emptyset.$ But this is a contradiction with
Lemma \ref{tan1}. 

Finally, in case (ii), we will show that
$W^s(\gamma)$ and $W^u(\eta)$ do not intersect in a transversal way. Let
$x\in W^s(\gamma)\cap W^u(\eta)$. We know that
$$T_xW^s(\gamma)=\langle X(x) \rangle\oplus E^1_x,\ \mbox{and}\
T_xW^u(\eta)=\langle X(x) \rangle\oplus E_x^2.$$ Hence, ${\rm
dim}(T_xW^s(\gamma)+T_xW^u(\eta))<{\rm dim}W^s(\gamma)+{\rm
dim}W^u(\eta)={\rm dim}M.$ This means that $W^s(\gamma)$ and
$W^u(\eta)$ do not intersect transversally. Since
$X\in\mathscr{K}$ we have $W^s(\gamma)\cap
W^u(\eta)=\emptyset.$ This is again a contradiction with Lemma
\ref{tan1}. $\square$\\

The following lemma is proved in \cite[Lemma 5.1]{ASS}.
\begin{lem} \label{differ} There is a residual set
$\mathscr{R}_3\subset\mathfrak{X}_{\mu}^1(M)$ (where $\dim M\geq 4$) such that if
$X\in\mathscr{R}_3,$ for any $C^1$-neighborhood $\mathcal{U}(X)$
of $X$ there are $Y\in\mathcal{U}(X)$ and two periodic orbits
$\gamma_Y, \eta_Y\in P_h(Y)$ such that ${\rm
index}(\gamma_Y)\not={\rm index}(\eta_Y),$ then $X$ has two
periodic orbits $\gamma, \eta\in P_h(X)$ with ${\rm
index}(\gamma)\not={\rm index}(\eta).$
\end{lem}

Note that if $p\in\gamma\in P(X)$, then $P_X^{\pi(p)}(p)$ has $1$
as eigenvalue with eigenvector $X(p),$ and all the other
eigenvalues are called the characteristic multipliers of $p.$

 For any $\delta>0$, we say that a point $p\in\gamma\in
P_h(X)$ has a {\it $\delta$-weak eigenvalue} if there is a
characteristic multiplier $\lambda$ of the orbit of $p$ associated to $P_X^t(p)$ such that
$(1-\delta)<|\lambda|<(1+\delta).$

\begin{lem}\label{nonhyper} Given $X\in \mathfrak{X}_{\mu}^1(M)$, if $\gamma\in P(X)$ is not hyperbolic, then there is $Y\in \mathfrak{X}_{\mu}^1(M)$,
$C^1$-close to $X$, such that $Y$ has orbits $\gamma_1, \gamma_2\in P_h(Y)$ with ${\rm
index}(\gamma_1)\not={\rm index}(\gamma_2)$.
\end{lem}

\noindent
 {\bf Proof.}  The proof is analog to the one in Lemma~\ref{PF} we just have to use Franks's lemma for conservative flows (\cite[Lemma 3.2]{BR0}). We leave the details for the reader. $\square$\\

\begin{lem}\label{nodelta} There is a residual set
$\mathscr{R}_4\subset\mathfrak{X}_{\mu}^1(M)$ such that for any
$X\in\mathscr{R}_4,$ if $X$ has the shadowing property, then there
is $\delta>0$ such that every $p\in\gamma\in P_h(X)$ does not have
a $\delta$-weak eigenvalue.
\end{lem}

\noindent
 {\bf Proof.} Let $\mathscr{R}_4=\mathscr{R}_2\cap\mathscr{R}_3,$
 and let $X\in\mathscr{R}_4$ satisfies the shadowing property. Suppose
 that there is $p\in\gamma \in P(X)$ such that for any $\delta>0$,
 $p$ has a $\delta$-weak eigenvalue. By  Lemma \ref{nonhyper}, we can
 find $Y$ $C^1$-close to $X$ such that $Y$ has two hyperbolic
 periodic orbits $\gamma_1, \gamma_2\in P(Y)$ with ${\rm
 index}(\gamma_1)\not={\rm index}(\gamma_2)$. But $X\in\mathscr{R}_2$ and  by Lemma \ref{same} this is a contradiction.
 $\square$\\

\begin{rk}\label{expansive2} 
Previous result also holds if we replace \emph{shadowing property} by \emph{expansive property}. Actually the proof in ~\cite[Lemma 2.6]{Se} can be followed step by step and we just have to use the Franks' lemma for conservative flows proved in ~\cite[Lemma 3.2]{BR0}.
\end{rk}

 A vector field $X\in\mathfrak{X}^1_{\mu}(M)$ is {\it a
divergence-free star vector field} if there exist a
$C^1$-neighborhood $\mathcal{U}(X)$ of $X$ such that for any
$Y\in\mathcal{U}(X)$, any element of the set $Crit(X)$, i.e., singularity or closed orbit, is
of hyperbolic type. The set of star divergence-vector fields is denoted by
$\mathcal{G}_{\mu}^\star(M)$.

The following result proved by Ferreira, is the multidimensional version of \cite{BR1}, and will be crucial in order to prove Theorem~\ref{thm2}.

\begin{lem}{\em \cite[Theorem 1]{F}}\label{star1} If
$X\in\mathcal{G}^\star_{\mu}(M)$, then $\text{Sing}(X)=\emptyset$ and $X$ is
Anosov.
\end{lem}

Now, we obtain the flow version of Lemma~\ref{two} whose proof follows the same steps of the discrete one and we leave it to the reader (see also \cite[Lemma 5.3]{ASS} and ~\cite[Lemma 5.1]{A}).

\begin{lem}\label{2delta} There is a residual set
$\mathscr{R}_5\subset\mathfrak{X}_{\mu}^1(M)$ such that for any
$X\in\mathscr{R}_5,$ if for any $C^1$-neighborhood
$\mathcal{U}(X)$ of $X$, there exist $Y\in\mathcal{U}(X)$ and
$p_Y\in\gamma_Y\in P_h(Y)$ such that $p_Y$  has a $\delta$-weak
eigenvalue, then $p\in\gamma\in P_h(X)$ has a $2\delta$-weak
eigenvalue.
\end{lem}

\noindent {\bf Proof of Theorem \ref{thm2}.} The proof is by contradiction. Let $\mathscr{R}_6=\mathscr{R}_4\cap
\mathscr{R}_5,$ and let $X\in\mathscr{R}_6$ and displaying also the
shadowing property. Suppose that $X\not\in\mathcal{G}^\star_{\mu}(M).$
Then there is a periodic point $p\in\gamma\in P(X)$
such that $p$ is has a $\delta/2$-weak eigenvalue. Since
$X\in\mathscr{R}_6,$ $p$ has a $\delta$-weak eigenvalue. This contradicts Lemma \ref{nodelta}. $\square$\\

\noindent {\bf Proof of Theorem \ref{thm22}.} Once again the proof is by contradiction. The residual subset stated in the theorem is defined by the intersection of $\mathscr{R}_5$ with the residual given by Remark~\ref{expansive2}. We let $X$ be expansive and in this residual. If $X\not\in\mathcal{G}^\star_{\mu}(M)$, then there is a periodic point $p\in\gamma\in P(X)$
such that $p$ is has a $\delta/2$-weak eigenvalue, but this contradicts Remark~\ref{expansive2}. $\square$\\

Now we consider the flow version of Lemma~\ref{manifolds} (see also \cite[Lemma 2.1]{SSY} and ~\cite[Theorem 3.3]{ASS}).

\begin{lem}\label{manifolds2}
Let $X\in \mathfrak{X}_{\mu}^1(M)$. If $\gamma_1,\gamma_2\in P_h(X)$ and $X$ has the weak specification property, then $W^u(\gamma_1,X)\cap W^s(\gamma_2,X)\not=\emptyset$.
\end{lem}

In an analogous way we obtain:

\begin{lem}\label{spec2}
There is a residual subset of $\mathfrak{X}_{\mu}^1(M)$ such that, if $X$ is in this residual and has the weak shadowing property, then there
is $\delta>0$ such that for any $p\in P(X),$ $p$ does not have a
$\delta$-weak eigenvalue.
\end{lem}

\noindent
 {\bf Proof of Theorem \ref{thm33}.} We define the residual subset of the theorem by the intersection of $\mathscr{R}_5$ with the residual given by Lemma~\ref{spec2}. Let $X$ be a vector field with the weak specification property and in the previous defined residual. To derive a contradiction, we may assume that
 $X\not\in\mathcal{G}_{\mu}^\star(M).$ Then there is a periodic point $p$ such
 that $p$ has a
 $\delta/2$-weak eigenvalue.  By Lemma \ref{2delta}, $p$ has a $\delta$-weak
 eigenvalue which contradicts Lemma~\ref{spec2} and Theorem~\ref{thm33} is proved. $\square$\\

\section{$C^1$-stably expansiveness and $C^1$-stably weak specification}\label{stable}

Some $C^1$-stability properties for volume-preserving maps were treated briefly as a remark in \cite[\S 3]{ST}. Here we developed the full detailed proofs for the expansiveness and the weak specification properties.

\subsection{$C^1$-stably expansiveness}

We say that $f\in \text{Diff}_\mu^{\,\,1}(M)$ is \emph{$C^1$-stably expansive} if any  $g\in \text{Diff}_\mu^{\,\,1}(M)$ in a small $C^1$-neighborhood of $f$  is expansive. We begin by proving the following result:

\begin{maintheorem}\label{T1exp}
If $f\in \text{Diff}_\mu^{\,\,1}(M)$, then $f$ is a volume-preserving Anosov map if and only if $f$ is $C^1$-stable expansive.
\end{maintheorem}

\noindent {\bf Proof.} 
The ``only if" part follows directly from stable manifold theory. Let us prove the ``if" statement; given $f$ in the $C^1$-interior subset of expansive diffeomorphisms in $\text{Diff}_\mu^{\,\,1}(M)$ which we denote by $\mathcal{E}^1_\mu(M)$, we prove that it is Anosov. By Theorem \ref{Ar} it is sufficient to show that if $f\in \mathcal{E}^1_\mu(M)$, then $f\in \mathcal{F}^\star_\mu(M)$. So let us prove next that all the periodic orbits of $f$ are hyperbolic. By contradiction, let us assume that there exists $f\in \mathcal{E}^1_\mu(M)$ having a  non-hyperbolic orbit $p\in P(f)$ of period $\pi$, i.e., $Df^{\pi}(p)$ has an eigenvalue $\sigma$ such that $|\sigma|=1$.

First of all we need to $C^1$-approximate the volume-preserving diffeomorphism $f$ by a new one, $f_1$, which, in the local coordinates given by Moser's theorem (\cite{Mo}), is \emph{linear} in a neighborhood of the periodic orbit $p$. To perform this task, in the volume-preserving setting, and taking into account ~\cite[Theorem 3.6.]{AM}, it is required higher differentiability. 

Thus, if $f$ is of class $C^\infty$, take $g=f$, otherwise we use \cite{Av} in order to obtain a $C^\infty$ $h\in\mathcal{E}^1_\mu(M)$, arbitrarily $C^1$-close to $f$, and such that $h$ has a periodic orbit $q$, close to $p$, with period $\pi$. We observe that $q$ may not be the analytic continuation of $p$ and this is precisely the case when $1$ is an eigenvalue of the tangent map $Df^{\pi}(p)$. 

If $q$ is not hyperbolic take $g=h$. If $q$ is hyperbolic for $Dh^{\pi}(q)$, then, since $h$ is $C^1$-arbitrarily close to $f$, some eigenvalue of $Dh^{\pi}(q)$ has norm arbitrarily close to one. Thus we apply the Franks' lemma for volume-preserving diffeomorphisms (\cite[Proposition 7.4]{BDP}) to obtain a $g\in\mathcal{E}^1_\mu(M)$ of class $C^\infty$, $C^1$-close to $h$ and such that $q$ is  a non-hyperbolic periodic orbit.

Now, we use the \emph{weak pasting lemma} (\cite[Theorem 3.6]{AM})  in order to obtain $f_1\in \mathcal{E}^1_\mu(M)$  such that, in local canonical coordinates, $f_1$ is linear and equal to $Dg$ in a neighborhood of the periodic non-hyperbolic orbit, $q$. Moreover, the existence of an eigenvalue, $\sigma$, with modulus equal to one is associated to an invariant subspace contained in the subspace $E^c_q\subseteq T_q M$ associated to norm-one eigenvalues. Furthermore, up to a perturbation using again ~\cite[Proposition 7.4]{BDP}, $\sigma$ can be taken rational. This fact assures the existence of periodic orbits arbitrarily close to the $f_1$-orbit of $q$. Thus, there exists $m\in\mathbb{N}$ such that $f_1^{m\pi}(q)|_{E^c_q}=(Dg^{m\pi})_q|_{E^c_q}=id$ holds, say in a $\xi$-neighborhood of $q$.

Finally, we just have to pick two points $x,y\in E^c_q$ sufficiently close in order to obtain $d(g^n(x),g^n(y))<\epsilon$ for all $n\in\mathbb{Z}$. It is clear that $g$ cannot be expansive which is a contradiction and Theorem~\ref{T1exp} is proved.
 $\square$\\

\begin{maintheorem}\label{T2exp}
If $X\in \mathfrak{X}_{\mu}^1(M)$, then $X$ is a divergence-free Anosov vector field if and only if $X$ is $C^1$-stable expansive.
\end{maintheorem}

\noindent {\bf Proof.} 
The proof is quite similar to the previous one. We just have to switch to the perturbation framework developed in \cite{AM} and ~\cite{BR0}. Namelly, the smoothness result ~\cite[Theorem 2.2]{AM}, the $C^{1+}$-pasting lemma for vector fields ~\cite[Theorem 3.1]{AM} and the Franks' lemma for divergence-free vector fields proved in ~\cite[Lemma 3.2]{BR0}. We leave to the reader the filling up all the details. $\square$\\

\subsection{$C^1$-stably weak specification}
We say that  $f\in \text{Diff}_\mu^{\,\,1}(M)$ has the \emph{$C^1$-stable weak specification property} if there exists a neighborhood $U(f)$ of $f$ in  $\text{Diff}_\mu^{\,\,1}(M)$  such that any $g\in U(f)$ satisfy the weak specification property\footnote{\,The results presented here can be stated for basic sets as it is done in \cite{SSY}. The changes do not involve extra difficulties and are straightforward.}. 

\begin{maintheorem}\label{T1spe}
If $f\in  \text{Diff}_\mu^{\,\,1}(M)$, then $f$ is a volume-preserving Anosov diffeomorphism if and only if $f$ has the $C^1$-stable weak specification property.
\end{maintheorem}

\subsubsection{Surface case} We begin by proving previous result for surfaces by using a symplectic reasoning. We prove that, if $f$ is an area-preserving map exhibiting the $C^1$-stable weak specification property, then $f$ is Anosov. If $f$ has the $C^1$-stable weak specification property and it is not Anosov, then, by ~\cite{N}, it can be $C^1$-approximated by an area-preserving map $f_0$ with elliptic orbits. By, ~\cite{Z} we can $C^1$-approximate $f_0$ by an area-preserving $C^\infty$ map $f_1$ with elliptic orbits and still with $C^1$-stable weak specification. Then, $f_1$ is $C^1$-approximated by area-preserving maps displaying invariant circles associated to irrational rotations and cannot by topologically mixing, thus, by ~\cite[Proposition 21.3]{DGS} cannot have the weak specification property.

\subsubsection{General case} Now, we consider manifolds with dimension $d\geq 3$.

In the next lemma we borrow the arguments in \cite[Lemma 2.2]{SSY} and combine it with a quite often result used along this paper, namely the Robinson's volume-preserving version of the Kupka-Smale theorem.

\begin{lem}\label{index}
If $f$ is a volume-preserving diffeomorphism having the $C^1$-stable weak specification property, then any two hyperbolic periodic points of $g\in U(f)$ have the same index.
\end{lem}
\noindent {\bf Proof.}  
Let $f\in\text{Diff}^{\,\,1}_\mu(M)$ having the $C^1$-stable weak specification property and $g\in U(f)$. Take any two hyperbolic periodic points $p$ and $q$ of $g$. Consider a $C^1$-neighborhood of $g$, $U_g$ such that any $h\in U_g$ has a well-defined analytic continuation of $p$ and $q$, denoted by $p_h$ and $q_h$ respectively.

Assuming that ${\rm index}(p)<{\rm index}(q)$ we will obtain a contradiction. By our assumption we obtain that $\dim W^u(p,g)+\dim W^s(q,g)<d$ where we recall that $d=\dim M$. Since, by Robinson's volume-preserving version of the Kupka-Smale theorem the stable/unstable manifolds intersects in a transversal way, any Robinson's (Kupka-Smale) volume-preserving diffeomorphism $h\in U_g$ must satisfy the condition 
\begin{equation}\label{W}
W^u(p_h,h)\cap W^s(q_h,h)=\emptyset,
\end{equation}
otherwise the intersection would not be transversal. But Lemma~\ref{manifolds} contradicts (\ref{W}).
$\square$\\

Since the ``only if'' part of Theorem~\ref{T1spe} is immediate, the proof of the ``if'' part is a direct consequence of the Theorem~\ref{Ar} and the following simple result:

\begin{lem}\label{main} 
If a volume-preserving diffeomorphism $f$ has the $C^1$-stably weak specification pro\-per\-ty, then $f\in\mathcal{F}^\star_{\mu}(M)$.  
\end{lem}

\noindent {\bf Proof.}  
The proof is by contradiction; let us assume that there exists a volume-preserving diffeomorphism $f$ having the $C^1$-stably specification property on a set $U(f)$ and exhibits a non-hyperbolic closed orbit $p$ of period $\pi$.

Using Lemma~\ref{PF} there exists $g\in U(f) $  such that $g$ has two hyperbolic periodic orbits with different indices which contradicts Lemma~\ref{index}.
$\square$\\

\begin{maintheorem}\label{T2spe}(\cite[Theorem 9.1]{ASS})
If $X\in \mathfrak{X}_{\mu}^1(M)$, then $X$ is a volume-preserving Anosov flow if and only if $X$ is has the $C^1$-stable weak specification property.
\end{maintheorem}

\noindent {\bf Proof.}  
We know that $X^t$ has the weak specification property, then $X^t$ is topologically mixing. Hence, $X^t$ has the $C^1$-stable weak specification property, then $X^t$ is $C^1$-robustly topologically mixing. By the result in ~\cite{BR0} $M$ displays a dominated splitting\footnote{If $d=3$, then $X^t$ is Anosov by ~\cite{AM}.}. Finally, an analogous argument for the absence of different indices allow us to conclude the global hyperbolicity.
$\square$\\

We observe that the proof of Theorem~\ref{T2spe} is different from the one followed in the diffeomorphism case. However, we could opt to use the same scheme by using the result in ~\cite[Theorem 2.2]{AM} and ~\cite{BD}.

\bigskip
\noindent {\bf Acknowledgement.} We would like to thank Alexander Arbieto by clarifying some aspects on his paper ~\cite{ASS} and also to Joana Torres for the careful reading of the manuscript. MB was partially supported by National Funds through FCT - ``Funda\c{c}\~{a}o para a Ci\^{e}ncia e a Tecnologia'', project PEst-OE/MAT/UI0212/2011. ML was supported by Basic
Science Research Program through the National Research Foundation
of Korea (NRF) funded by the Ministry of Education, Science and
Technology, Korea (No. 2011-0007649).

\end{document}